\documentclass[a4paper, 11pt]{amsart}
\usepackage{amssymb,amsmath,txfonts,mathrsfs,mathtools,titletoc}
\usepackage[alphabetic]{amsrefs}
\newtheorem{theorem}{Theorem}[section]
\newtheorem{prop}[theorem]{Proposition}
\newtheorem{lemma}[theorem]{Lemma}
\newtheorem{remark}[theorem]{Remark}
\newtheorem{claim}[theorem]{Claim}
\newtheorem{question}[theorem]{Question}
\newtheorem{definition}[theorem]{Definition}
\newtheorem{cor}[theorem]{Corollary}
\newtheorem{example}[theorem]{Example}

\numberwithin{equation}{section}
\def\pf{{\it Proof:}~}
 
\begin{document}

\title[Lower bound of Ricci flow's Existence time]{Lower bound of Ricci flow's Existence time}
\author{Guoyi Xu}
 \address{Yau Mathematical Sciences Center \\ Tsinghua University, Beijing\\P. R. China}
\email{gyxu@math.tsinghua.edu.cn}
\date{\today}

\begin{abstract}
Let $(M^n, g)$ be a compact $n$-dim ($n\geq 2$) manifold with nonnegative Ricci curvature, and if $n\geq 3$ we assume that $(M^n, g)\times \mathbb{R}$ has nonnegative isotropic curvature. The lower bound of the Ricci flow's existence time on $(M^n, g)$ is proved. This provides an alternative proof for the uniform lower bound of a family of closed Ricci flows' maximal existence times, which was firstly proved by E. Cabezas-Rivas and B. Wilking. We also get an interior curvature estimates for $n= 3$ under $Rc\geq 0$ assumption among others. Combining these results, we proved the short time existence of the Ricci flow on a large class of $3$-dim open manifolds, which admit some suitable exhaustion covering and have nonnegative Ricci curvature.
\\[3mm]
\end{abstract}

\thanks{The author was partially supported by NSFC 11401336}

\maketitle

\section{Introduction}
In \cite{Ham}, R. Hamilton introduced Ricci flow and proved that the Ricci flow on any compact $3$-dimensional manifold $(M^3, g)$ with $Rc> 0$ will deform the metric $g$ to a metric with constant positive sectional curvature. In \cite{Ham88}, he further proved the above theorem for $2$-dim compact surfaces with positive curvature. When the dimension of manifold $n\geq 4$, the isotropic curvature introduced by M. Micallef and J. Moore \cite{MM} is an important concept. If $(M^n, g)\times \mathbb{R}$ has positive isotropic curvature, S. Brendle \cite{Brendle} proved that the Ricci flow deforms $(M^n, g)$ to a Riemannian manifold of constant positive sectional curvature. 

From the definition of isotropic curvature, when $n= 3$, that $(M^n, g)\times \mathbb{R}$ has nonnegative isotropic curvature is equivalent to $Rc\geq 0$ (also see \cite{Brendle}). From all the above, if $Rc> 0$ when $n= 2$ and $(M^n, g)\times \mathbb{R}$ has positive isotropic curvature when $n\geq 3$, Ricci flow deforms $(M^n, g)$ to a Riemannian manifold of constant positive sectional curvature.

In all the above case, using maximal principle, it is easy to get the upper bound of existence time of the Ricci flow on $(M^n, g)$. In \cite{Simon}, M. Simon posed the following question (Problem $3.1$ there):
\begin{question}\label{ques 1}
What elements of the geometry need to be controlled, in order to guarantee
that a solution of Ricci flow does not become singular?
\end{question}

In this note, we give a lower bound of existence time of the Ricci flow on compact manifold $(M^n, g)$ with $Rc\geq 0$, and if $n\geq 3$ we assume that $(M^n, g)\times \mathbb{R}$ has nonnegative isotropic curvature. To state our result, we use $B_t(x, r)$ to denote the geodesic ball centered at $x$ with radius $r$ on $(M^n, g(t))$, and $V_{t}(B_s(x, r))$ is the volume of $B_s(x, r)$ with respect to the metric $g(t)$. One of our main results is the following:
\begin{theorem}\label{thm 1.1}
{Let $(M^n, g)$ be a compact $n$-dim ($n\geq 2$) manifold with $Rc\geq 0$, and if $n\geq 3$ we assume that $(M^n, g)\times \mathbb{R}$ has nonnegative isotropic curvature. There exists positive constants $\delta= \delta(n)$, $\epsilon_0= \epsilon_0(n)$, such that if $x_0\in M^n$ is a fixed point satisfying $|Rm|\leq r_0^{-2}$ on $B_{0}(x_0, r_0)$ and $B_{0}(x_0, r_0)$ is $\delta$-almost isoperimetrically Euclidean, where $r_0> 0$ is some constant; then the Ricci flow exists on $M^n \times [0, T]$ where $T= T(n, r_0)> 0$ and
\begin{equation}\nonumber
{|Rm|(x, t)\leq 2r_0^{-2}\ , \quad \quad x\in B_{t}\Big(x_0, \frac{1}{2}\epsilon_0 r_0\Big)\ , \ t\in [0, T] 
}
\end{equation}  
}
\end{theorem}

\begin{remark}\label{rem 1.1}
{One interesting thing of such lower bound of existence time is that it only depends on the local property of manifolds in a neighborhood of a point, not globally on the whole manifold. This kind of lower bound on existence time for Ricci flow had been got in Corollary $5$ of \cite{CW} for any complete manifolds with stronger curvature assumptions. Note for any point $x_0\in M^n$, we can always find the corresponding $r_0$ satisfying the assumptions in the above theorem, although such $r_0$ may be very small. This also provides a partial answer to Question \ref{ques 1} when manifolds satisfy the assumptions in the above theorem.
}
\end{remark}

The existence time's lower bound part of the following theorem is not totally new (see Corollary $5$ in \cite{CW}), however in $3$-dimension the results seems to be new. And our method to get such lower bound is different from the method 
in \cite{CW}.
\begin{theorem}\label{thm 1.2}
{Let $(M^n, g)$ be a compact $n$-dim manifold with $Rc\geq 0$ where $n\geq 2$. If $n\geq 4$, we assume $K_g^{\mathbb{C}}\geq 0$. Let $x_0\in M^n$ be a fixed point satisfying $|Rm|\leq r_0^{-2}$ on $B_{0}(x_0, r_0)$ and $V_{0}(B_{0}(x_0, r_0))\geq \nu_0r_0^n$ for some $\nu_0> 0$, $r_0> 0$. Then there exists $T= T(n, \nu_0, r_0)> 0$ such that the Ricci flow exists on $M^n \times [0, T]$. Also on $[0, T]$,
\begin{equation}\nonumber
{|Rm|(x, t)\leq 2r_0^{-2} \quad x\in B_{t}\Big(x_0, \frac{r_0}{2}\Big)
}
\end{equation}
Furthermore for any $r> \frac{r_0}{2}$, there is a $B= B(n, \nu_0, r_0, r)> 0$, such that on $[0, T]$, 
\begin{equation}\label{1.2.1}
{|Rm(x, t)|\leq B+ \frac{B}{t} \quad x\in B_{t}\big(x_0, \frac{r}{4}\big)
}
\end{equation}
}
\end{theorem}

\begin{remark}\label{rem 1.2}
{Note that $(M^n, g)\times \mathbb{R}$ has nonnegative isotropic curvature is implied by nonnegative complex sectional curvature when $n\geq 4$. Theorem \ref{thm 1.2} has stronger assumptions on curvature than Theorem \ref{thm 1.1}. The stronger assumptions about initial curvature is to guarantee getting (\ref{1.2.1}), which is very important for our application to short-time existence of Ricci flow on noncompact manifolds (see the proof of Corollary \ref{cor 4.2}).
}
\end{remark}

The local curvature estimates in Theorem \ref{thm 1.2} were known in $3$-dimension to B.-L. Chen (Theorem $3.6$ of \cite{Chen}). Under different assumptions, the similar local curvature estimates for $2$-dim and $3$-dim Ricci flows were obtained by M. Simon (Theorems $1.1$ and $1.5$ in \cite{Simon3}). In \cite{Simon2} (Lemma $4.3$ and Theorem $7.1$ there), he obtained the same result as (\ref{1.2.1}), provided that Ricci curvature has lower bound, and every geodesic ball of radius $r_0$ has a uniform bound from below independent of the center of the ball. 

One motivation to get such lower bound of existence time is for studying the short-time existence of Ricci flow on open manifold $(M^n, g)$ with the assumption in the above theorems.

When $(M^n, g)$ has nonnegative complex sectional curvature, by considering the doubling of convex sets contained in a Cheeger-Gromoll convex exhaustion of $(M^n, g)$ and solving the initial value problem for the Ricci flow on these closed manifolds, E. Cabezas-Rivas and B. Wilking \cite{CW} obtain a sequence of closed solutions of the Ricci flow $(M_i, g_i(t))$ with nonnegative complex sectional curvature which subconverge to a solution of the Ricci flow $(M^n, g(t))$. Finally they successfully established the short-time existence of the Ricci flow on $(M^n, g)$ with $K_g^{\mathbb{C}}\geq 0$ among other things (Theorem $1$ in \cite{CW}).

One of the keys to get the above short-time existence result is to establish the following claim:
\begin{claim}[Proposition $4.3$ in \cite{CW}]\label{claim 1.1}
{There is a lower bound (independent of $i$) for the maximal times of existence $T_i$ for each Ricci flow $(M_i, g_i(t))$ obtained above.} 
\end{claim}

In \cite{CW}, to get the above lower bound, the authors used an upper bound of the scalar curvature integral on any unit geodesic ball of manifold with curvature bounded from below, which was established by Petrunin in Theorem $1.1$ of \cite{Petr}. Petrunin's beautiful result was proved by using deep results in Alexandrov spaces, although the conclusion is in the category of Riemannian geometry. 

One of the applications of Theorem \ref{thm 1.1} is to provide a direct proof of Claim \ref{claim 1.1} without involving the results in Alexandrov spaces. Moreover using Theorem \ref{thm 1.2} in $3$-dimensional case, and following the same strategy in \cite{CW}, we can get the short-time existence of Ricci flow on a large class of $3$-manifolds including all $3$-manifolds with nonnegative sectional curvature (note the nonnegative sectional curvature case was firstly proved in \cite{CW}). M. Simon (Theorem $1.9$ of \cite{Simon2}) proved the short-time existence of Ricci flow on a class of $3$-manifolds with $Rc\geq 0$ and satisfying other assumptions about curvature near infinity, the injective radius. It seems that our class of manifolds do not have much overlap with his class. 

A long term goal in $3$-dim case is to prove the short-time existence of Ricci flow on complete manifolds with $Rc\geq 0$. We hope that our result will be helpful on studying this long term goal. The more ambitious target is to study the short time existence of Ricci flow on any complete manifold $(M^n, g)$, provided that $(M^n, g)\times \mathbb{R}$ has nonnegative isotropic curvature. In this more general scheme, $3$-dim case will be exactly the $Rc\geq 0$ case.

\section{Perelman's pseudolocality and local Curvature estimates}

In the following computations,  we will use some cut-off functions which are composition of  a  cut-off function of $\mathbb R$ and distance function. A cut-off function $\varphi$ on real line $\mathbb R$, is a smooth nonnegative nonincreasing function, it is $1$ on $(-\infty, 1]$ and $0$ on $[2, \infty)$. We can further assume that
\begin{equation}\label{2.1.1}
{|\varphi '|\leq 2 \ , \quad |\varphi ''|+ \frac{(\varphi ')^2}{\varphi}\leq 16.
}
\end{equation}
Another often used notation is $\square= \frac{\partial}{\partial t}- \Delta$, where $\Delta$ is the Laplacian with the metric $g(t)$.

Let us recall Perelman's pseudolocality theorem. We have the following definition:
\begin{definition}\label{def almost isop}
{We say that $B_{g(0)}(x_0, r_0)$ is \textbf{$\delta$-almost isoperimetrically Euclidean} if
\begin{equation}\nonumber
{\Big( V_0(\partial \Omega) \Big)^n\geq (1- \delta)c_n\Big(V_0(\Omega)\Big)^{n- 1}
}
\end{equation}
for any regular domain $\Omega\subset B_0(x_0, r_0)$, where $c_n\vcentcolon= n^n\omega_n$ is the Euclidean isoperimetric constant and $\omega_n$ is the volume of the Euclidean unit $n$-ball. 
}
\end{definition}

\begin{theorem}[Perelman's pseudolocality, \cite{Pere}]\label{thm pseu}
{For every $\alpha> 0$ and $n\geq 2$, there exist $\delta= \delta(\alpha, n)> 0$ and $\epsilon_0= \epsilon(\alpha, n)> 0$ with the following property. If $(M^n, g(t))_{t\in [0, T]}$ is a smooth solution of the Ricci flow on compact manifold $(M^n, g(0))$, $R(x, 0)\geq -r_0^{-2}$ for all $x\in B_{g(0)}(x_0, r_0)$ where $r_0> 0$ and $x_0\in M^n$ are fixed, and $B_{g(0)}(x_0, r_0)$ is $\delta$-almost isoperimetrically Euclidean; then we have
\begin{equation}\nonumber
{|Rm|(x, t)\leq \frac{\alpha}{t}+ \frac{1}{\big(\epsilon_0 r_0\big)^2}\ , \quad \quad  x\in B_t(x_0, \epsilon_0 r_0)\ , \ t\in \Big(0, \min\{\big(\epsilon_0 r_0\big)^2, T\}\Big]
}
\end{equation}
}
\end{theorem}

\begin{remark}\label{rem per}
{If we choose $\alpha= 2$ in the above theorem, fix $n$, $x_0$ and suitable $r_0$, we get that there exists numerical constants $\delta= \delta(2, n)$ and $\epsilon_0= \epsilon_0(2, n)< 1$ such that 
\begin{equation}\label{2.1}
{|Rm|(x, t)\leq \frac{2}{t}+ \frac{1}{\big(\epsilon_0 r_0\big)^2}\ , \quad \quad  x\in B_t(x_0, \epsilon_0 r_0)\ , \ t\in \Big(0, \min\{\big(\epsilon_0 r_0\big)^2, T\}\Big]
}
\end{equation}
}
\end{remark}

We also need a theorem due to B.-L. Chen (see Theorem $3.1$ in \cite{Chen}), and similar results were also obtained independently by M. Simon (Lemma $6.2.2$ in \cite{SimonThesis}, Theorem $1.3$ in \cite{SimonIMRN}). 
\begin{theorem}[B. -L.  Chen, M. Simon]\label{thm chen}
{There is a constant $C= C(n)$ with the following property. Suppose we have a smooth Ricci flow solution $(M^n, g(t))_{t\in [0, T]}$ such that $B_{t}(x_0, r_0)$, $0\leq t\leq T$, is compactly contained in $M^n$, with
\begin{equation}\nonumber
{|Rm|(x, 0)\leq r_0^{-2}\ , \quad \quad \quad  \ x\in B_{0}(x_0, r_0)
}
\end{equation}
and for some $K\geq 1$,
\begin{equation}\nonumber
{|Rm|(x, t)\leq \frac{K}{t}\ , \quad \quad \quad x\in B_t(x_0, r_0)\ , \quad t\in (0, T]
}
\end{equation}
Then we have 
\begin{equation}\nonumber
{|Rm|(x, t)\leq e^{CK}\Big(r_0- d_t(x_0, x)\Big)^{-2}\ , \quad \quad \quad x\in B_t(x_0, r_0)\ , \quad t\in [0, T]
}
\end{equation}
}
\end{theorem}

The following proposition is motivated by Theorem $3.6$ in \cite{Chen}. The strategy of the proof is applying Perelman's pseudolocality theorem to get (\ref{2.1}), then we can use Theorem \ref{thm chen} to get (\ref{2.2.1}).

\begin{prop}\label{prop 2.2}
{Let $(M^n, g(t))_{t\in [0, T]}$ be a smooth solution of the Ricci flow on compact manifold $(M^n, g(0))$. There exists positive constants $\delta= \delta(n)$, $\epsilon_0= \epsilon_0(n)$, such that if $x_0\in M^n$ is a fixed point satisfying $|Rm|\leq r_0^{-2}$ on $B_{0}(x_0, r_0)$ and $B_{0}(x_0, r_0)$ is $\delta$-almost isoperimetrically Euclidean, where $r_0> 0$ is some constant; then there exists $C= C(r_0, n)> 0$ such that 
\begin{equation}\label{2.2.1}
{|Rm|(x, t)\leq 2r_0^{-2}\ , \quad \quad \quad x\in B_{t}\Big(x_0, \frac{1}{2}\epsilon_0 r_0\Big)\ , \quad t\in \Big[0, \min\{T, C(r_0, n)\}\Big]
}
\end{equation}  
}
\end{prop}

\pf
{We define $T_0$ as the following:
\begin{align}
T_0\vcentcolon= \max{ \Big\{t\Big| \ t\leq T,\ |Rm|(x, s)\leq 2r_0^{-2},\ when \ x\in B_{s}\Big(x_0, \frac{1}{2}\epsilon_0 r_0\Big),\ and \ s\in [0, t]\Big\} } \nonumber 
\end{align}
Recall by assumption, $|Rm(x, 0)|\leq r_0^{-2}$ on $B_{0}(x_0, r_0)$. Without loss of generality, we assume that $T_0< \min\{T, 1\}$. Then there is $(x_1, t_1)$ such that $t_1\leq T_0$, $x_1\in B_{t_1}\Big(x_0, \frac{1}{2}\epsilon_0 r_0\Big)$, $|Rm|(x_1, t_1)= 2r_0^{-2}$. 
By Theorem \ref{thm pseu} and Remark \ref{rem per}, there exists $\delta= \delta(n)$, $\epsilon_0= \epsilon_0(n)$, such that for $x\in B_{t}(x_0, \epsilon_0 r_0)$, $t\in (0, \min\{T_0, \big(\epsilon_0 r_0\big)^2\}]$,
\begin{equation}\label{2.2}
{|Rm(x, t)|\leq \frac{2+ \frac{T_0}{\big(\epsilon_0 r_0\big)^2}}{t}\leq \frac{C(n, r_0)}{t}
}
\end{equation}

Now apply Theorem \ref{thm chen}, combining with (\ref{2.2}), we have $|Rm|(x, t)\leq C(n, r_0)$ on $B_t\Big(x_0, \frac{3}{4}\epsilon_0 r_0\Big)$, where $t\in \Big[0, \min\{T_0, \big(\epsilon_0 r_0\big)^2\}\Big]$.  
 
Consider the evolution equation of $u= \varphi(\frac{4d_t(x_0, x)}{\epsilon_0 r_0}- 1)|Rm(x, t)|^2$, where $\varphi$ is the smooth nonnegative  decreasing function chosen as in the beginning of this section. 
\begin{equation}\nonumber
{\left.
\begin{array}{rl}
\square u &\leq \Big(\frac{4\varphi '}{r_0}\square d_{t}(x_0, x)- \frac{16}{r_0^2}\varphi ''\Big)|Rm|^2- 2\varphi |\nabla Rm|^2 \\
& -2\nabla \varphi \nabla |Rm|^2 + C_1\varphi |Rm|^3
\end{array} \right.
}
\end{equation}
where $C_1$ is some constant depending only on $n$. 

Let $u_{max}(t)= u(x(t), t)= \max_{x\in M^3}u(x, t)$, then by the maximum principle and (\ref{2.1.1}), and using Lemma $8.3$ in \cite{Pere}, at $(x(t), t)$ we have
\begin{equation}\nonumber
{\left.
\begin{array}{rl}
\frac{d^{-}}{dt}u_{max}&\leq C(n, r_0)|Rm|^2+ C|Rm|u_{max}\\
&\leq C(n, r_0)+ C(n, r_0)u_{max}= C+ Cu_{max}
\end{array} \right.
}
\end{equation}
where $C= C(n, r_0)$ only depends on $n, r_0$. Integrating this inequality, 
\begin{equation}\nonumber
{u_{max}(t)\leq [u_{max}(0)+ 1]e^{Ct}- 1\leq [r_0^{-4}+ 1]e^{Ct}- 1
}
\end{equation}

Hence there exists $C(n, r_0)> 0$ such that 
\begin{align}
u_{max}(t)< 4r_0^{-4}\ , \quad  0\leq t\leq \min\{ C(n, r_0), T_0\} \nonumber
\end{align}
Then we get $|Rm(x, t)|< 2r_0^{-2}$ on $B_{t}\Big(x_0, \frac{1}{2}\epsilon_0 r_0\Big)$ when $0\leq t\leq \min\{C(\nu_0, r_0), T_0\}$. 
Recall that there is $(x_1, t_1)$ such that 
\begin{align}
t_1\leq T_0\ , \quad x_1\in B_{t_1}\Big(x_0, \frac{1}{2}\epsilon_0 r_0\Big)\ , \quad |Rm|(x_1, t_1)= 2r_0^{-2} \nonumber
\end{align}
Therefore $T_0\geq \min\{T, C(n, r_0) \}$. 
}
\qed

We recall the definition of complex sectional curvature.
\begin{definition}\label{def 2.1}
{Let $(M^n, g)$ be a Riemannian manifold, and consider its complexified
tangent bundle $T^{\mathbb{C}}M := TM\otimes \mathbb{C}$. We extend the curvature tensor $Rm$ and the metric $g$ at $p$ to $\mathbb{C}$-multilinear maps $Rm: (T^{\mathbb{C}}_p M)^4 \rightarrow \mathbb{C}$, $g: (T^{\mathbb{C}}
_p M)^2 \rightarrow \mathbb{C}$. The complex sectional curvature of a $2$-dimensional complex subspace $\sigma$ of $T^{\mathbb{C}}_p M$ is defined by
\begin{align}
K^{\mathbb{C}}(\sigma)= Rm(u, v, \bar{v}, \bar{u})= g\Big(Rm(u\wedge v), \overline{u\wedge v}\Big), \nonumber
\end{align}
where $u$ and $v$ form any unitary basis for $\sigma$, i.e. $g(u, \bar{u})= g(v, \bar{v})= 1$ and $g(u, \bar{v})= 0$. We say $(M^n, g)$ has nonnegative complex sectional curvature if $K_g^{\mathbb{C}}\geq 0$.
}
\end{definition}

The following proposition is very close to Proposition \ref{prop 2.2}, and the assumption about initial curvature is to guarantee getting (\ref{3.2.2}), which is important for our application to the Ricci flow on noncompact manifolds.
\begin{prop}\label{prop 3.2}
{Let $(M^n, g(x, 0))$ be a compact $n$-dim manifold with $Rc\geq 0$. If $n\geq 4$, assume $K_g^{\mathbb{C}}\geq 0$. Let $x_0\in M^n$ be a fixed point satisfying $|Rm|\leq r_0^{-2}$ on $B_{0}(x_0, r_0)$ and $V_{0}(B_{0}(x_0, r_0))\geq \nu_0r_0^n$ for some $\nu_0> 0$, $r_0> 0$. Let $g(x, t)$ , $t\in [0, T]$ be a smooth solution to the Ricci flow with $g(x, 0)$ as initial metric. Then there exists $C= C(n, \nu_0, r_0)> 0$ such that 
\begin{equation}\nonumber
{|Rm(x, t)|\leq 2r_0^{-2}\ , \quad \quad  x\in B_{t}\Big(x_0, \frac{r_0}{2}\Big) \ , \quad  t\in \big[0, \min\{T, C(n, \nu_0, r_0)\}\big] 
}
\end{equation}  
Furthermore, for any $r> \frac{r_0}{2}$, there exists $B= B(n, \nu_0, r_0, r)> 0$ such that 
\begin{equation}\label{3.2.2}
{|Rm(x, t)|\leq B+ \frac{B}{t}\ , \quad \quad  x\in B_{t}\Big(x_0, \frac{r}{4}\Big) \ , \quad t\in \big[0, \min\{T, C(n, \nu_0, r_0)\}\big] 
}
\end{equation}
}
\end{prop}

\begin{remark}\label{rem 2.2}
{This proposition generalizes the compact case of Theorem $3.6$ in \cite{Chen} to higher dimensions. The observation is: we need two key facts in the proof of Theorem $3.6$ there. The first one is $Rc\geq 0$ during Ricci flow. The other one is that asymptotic volume ratios of related ancient solutions with the same curvature assumptions as initial metric are $0$. When we can get these two facts from initial condition, our proof follows the lines given in the proof of Theorem $3.6$ in \cite{Chen}, although a number of modifications are necessary. Similar result as (\ref{3.2.2}) was previously obtained in Proposition $4.6$ of \cite{CW}, although our proof is somewhat different from theirs.
}
\end{remark}

\pf
{From the assumption, we get that $Rc\geq 0$ is preserved from \cite{Ham} when $n=2, 3$. When $n\geq 4$, $K_g^{\mathbb{C}}\geq 0$ is preserved from \cite{BS}. Note $K_g^{\mathbb{C}}\geq 0$ implies $Rc\geq 0$, hence $Rc(x, t)\geq 0$ on $M^n\times [0, T]$.

We define $T_0$ as the following:
\begin{align}
T_0\vcentcolon= \max{ \Big\{t\Big| \ t\leq T,\ |Rm|(x, s)\leq 2r_0^{-2},\ if \ x\in B_{s}\Big(x_0, \frac{r_0}{2}\Big)\ , \ s \in [0, t]\Big\} } \nonumber 
\end{align}
By assumption, $|Rm(x, 0)|\leq r_0^{-2}$ on $B_{0}(x_0, r_0)$. Without loss of generality, we assume that $T_0< T$. Then there is $(x_1, t_1)$ such that $t_1\leq T_0$, $x_1\in B_{t_1}(x_0, \frac{r_0}{2})$, $|Rm|(x_1, t_1)= 2r_0^{-2}$. 
We will estimate $T_0$ from below by a positive constant depending only on $n$,  $\nu_0$ and $r_0$. We have the following claim:
\begin{claim}\label{claim 3.3}
{For any $r> \frac{r_0}{2}$, there is a $B= B(n, \nu_0, r_0, r)> 0$, such that 
\begin{equation}\label{3.3.1}
{|Rm(x, t)|\leq B+ \frac{B}{t}\ , \quad \quad \quad x\in B_{t}\Big(x_0, \frac{r}{4}\Big)\ , \quad t\in [0, T_0]
}
\end{equation}

}
\end{claim}
{\it \textbf{Proof of Claim \ref{claim 3.3}}:}~
{We argue by contradiction. Suppose (\ref{3.3.1}) does not hold. 

Then for some fixed $r> \frac{r_0}{2}$, there exists a sequence of Ricci flow solutions $(M_i^n, g^{(i)}(t))$ on $[0, T_{0,i}]$, satisfying the assumptions (using $x_0^{(i)}$ to replace $x_0$ there) in Proposition \ref{prop 3.2} and
\begin{align}
|Rm(x, s)|_{g^{(i)}(t)}\leq 2r_0^{-2}\ , \quad \quad  x\in B_{g^{(i)}(s)}\Big(x_0^{(i)}, \frac{r_0}{2}\Big)\ , \quad s\in [0, T_{0,i}] \label{contradiction assumption 1}
\end{align}
Furthermore, there exists some $x_i\in B_{g^{(i)}(t_i)}\big(x_0^{(i)}, \frac{r}{4}\big)$, $t_i\in [0, T_{0,i}]$, such that
\begin{align}
|Rm(x_i, t_i)|_{g^{(i)}(t_i)}\geq B_i+ \frac{B_i}{t_i}\ , \quad \quad \quad \lim_{i\rightarrow \infty}B_i= \infty \nonumber
\end{align}

Using Lemma \ref{lem 3.1} in the Appendix, let $A_i= \frac{r}{25}B_i^{\frac{1}{2}}\geq 1$, we can choose $\bar{x}_i\in B_{g^{(i)}(\bar{t}_i)}\Big(x_0^{(i)}, 2A_iB_i^{-\frac{1}{2}}+ \frac{r}{4}\Big)$, $\bar{t}_i \in [0, T_{0,i}]$ with $\bar{Q}_i \vcentcolon = |Rm(\bar{x}_i, \bar{t}_i)|_{g^{(i)}(\bar{t}_i)}\geq B_i\bar{t}_i^{-1}$,
\begin{align} 
|Rm(x, t)|_{g^{(i)}(t)}\leq 4\bar{Q}_i \ , \quad (x, t)\in B_{g^{(i)}(\bar{t}_i)}\Big( \bar{x_i}, \frac{1}{10}A_i^{\frac{1}{2}}\bar{Q}_i^{-\frac{1}{2}}\Big) \times \Big[\bar{t}_i-\frac{2}{25n}A_i\bar{Q}_i^{-1}, \bar{t}_i\Big] \label{bound of Rm}
\end{align}

From the Ricci flow equation, 
\begin{align}
\frac{\partial}{\partial t}V_{g^{(i)}(t)}\Big(B_{g^{(i)}(0)}(x_0^{(i)}, \frac{r_0}{2})\Big)&= -\int_{B_{g^{(i)}(0)}(x_0^{(i)}, \frac{r_0}{2})} R(x, t)d\mu_{i}(x) \nonumber \\
&\geq -C(n)r_0^{-2}V_{g^{(i)}(t)}\Big(B_{g^{(i)}(0)}(x_0^{(i)}, \frac{r_0}{2})\Big) \nonumber
\end{align}
in the last inequality above we used (\ref{contradiction assumption 1}). Hence 
\begin{equation}\label{volume lower bound 0}
{V_{g^{(i)}(t)}\Big(B_{g^{(i)}(0)}(x_0^{(i)}, \frac{r_0}{2})\Big)\geq C\nu_{0} r_0^{n}e^{-Cr_0^{-2}}
}
\end{equation}
hereafter we denote $C\nu_{0} r_0^{n}e^{-Cr_0^{-2}}$ in (\ref{volume lower bound 0}) by $\phi(\nu_{0}, r_0)$.

From $Rc(x, t)\geq 0$ and the Ricci flow equation, we get 
\begin{align}
B_{g^{(i)}(0)}(x_0^{(i)}, \frac{r_0}{2})\subset B_{g^{(i)}(t)}(x_0^{(i)}, \frac{r_0}{2})\ , \quad \quad t\in[0, T_{0,i}] \label{subset}
\end{align}

By (\ref{subset}) and the assumption $r> \frac{r_0}{2}$, also use (\ref{volume lower bound 0}),
\begin{equation}\label{3.3.2}
{V_{g^{(i)}(t)}(B_{g^{(i)}(t)}(x_0^{(i)}, r))\geq V_{g^{(i)}(t)}\Big(B_{g^{(i)}(0)}(x_0^{(i)}, \frac{r_0}{2})\Big)\geq \phi(\nu_{0}, r_0)\ ,  \quad \quad t\in [0, T_{0,i}] 
}
\end{equation}

Scaling metric by $\bar{g}^{(i)}(t)= \bar{Q}_ig^{(i)}(\bar{t}_i+ t\bar{Q}_i^{-1})$ around $(\bar{x}_i, \bar{t}_i)$. Note 
\begin{equation}\nonumber
{d_{g^{(i)}(\bar{t}_i)}(\bar{x}_i, x_0^{(i)})\leq \frac{r}{4}+ 2A_iB_i^{-\frac{1}{2}}< \frac{r}{3}
}
\end{equation}
which implies 
\begin{align}
B_{g^{(i)}(\bar{t}_i)}\big(x_0^{(i)}, \frac{r}{3}\big)\subset B_{g^{(i)}(\bar{t}_i)}(\bar{x}_i, \frac{2}{3}r)  \label{3.3.3}
\end{align}

If $0< s\leq \frac{2}{3}r\bar{Q}_i^{\frac{1}{2}}$, from $Rc(g^{(i)}(t))\geq 0$ for $t\in [0, T]$, using the Bishop-Gromov Volume Comparison Theorem, 
\begin{align}
&\frac{V_{\bar{g}^{(i)}(0)}\Big(B_{\bar{g}^{(i)}(0)}\big(\bar{x}_i, s\big)\Big)}{s^n}= \frac{V_{g^{(i)}(\bar{t}_i)} \Big(B_{g^{(i)}(\bar{t}_i)}\big(\bar{x}_i, s\bar{Q}_i^{-\frac{1}{2}}\big)\Big)}{(s\bar{Q}_i^{-\frac{1}{2}})^n} \geq \frac{V_{g^{(i)}(\bar{t}_i)}\big(B_{g^{(i)}(\bar{t}_i)}(\bar{x}_i, \frac{2}{3}r)\big)}{(\frac{2}{3}r)^n}  \nonumber \\
&\quad \geq \frac{V_{g^{(i)}(\bar{t}_i)}\Big(B_{g^{(i)}(\bar{t}_i)}\big(x_0^{(i)}, \frac{r}{3}\big)\Big)}{\Big(\frac{2}{3}r\Big)^n}\geq  \frac{1}{\Big(\frac{2}{3}r\Big)^n}\cdot\frac{V_{g^{(i)}(\bar{t}_i)}\Big(B_{g^{(i)}(\bar{t}_i)}\big(x_0^{(i)}, r\big)\Big)}{3^n}
\geq \frac{\phi(\nu_0, r_0)}{(2r)^n} \label{3.3.4}
\end{align}
we used (\ref{3.3.3}) in the second inequality above, (\ref{3.3.2}) was used in the last inequality.

Note (\ref{3.3.4}) will imply
\begin{align}
V_{\bar{g}^{(i)}(0)}\Big(B_{\bar{g}^{(i)}(0)}\big(\bar{x}_i, s\big)\Big)\geq \Big(\frac{s}{2r}\Big)^n\phi(\nu_0, r_0)\ , \quad \quad s\in [0, \frac{2}{3}r\bar{Q}_i^{\frac{1}{2}}] \label{volume lower bound}
\end{align}

From (\ref{bound of Rm}) and the definition of $\bar{g}^{(i)}$, 
\begin{align}
|\overline{Rm}(x, t)|_{\bar{g}^{(i)}(t)}\leq 4 \ , \quad\quad  x\in B_{\bar{g}^{(i)}(t)}\Big(\bar{x}_i, \frac{1}{10}A_i^{\frac{1}{2}}\Big)\ , \quad t\in \Big[-\frac{2}{25n}A_i, 0\Big] \label{3.3.5} 
\end{align}

Let $t= 0$ in (\ref{3.3.5}), we get 
\begin{align}
|\overline{Rm}(x, 0)|_{\bar{g}^{(i)}(0)}\leq 4 \ , \quad\quad  x\in B_{\bar{g}^{(i)}(0)}\Big(\bar{x}_i, \frac{1}{10}A_i^{\frac{1}{2}}\Big) \label{3.3.7}
\end{align}

From (\ref{volume lower bound}) and (\ref{3.3.7}), using Theorem $4.3$ in \cite{CGT}, we get 
\begin{align}
inj_{\bar{g}^{(i)}(0)}(\bar{x}_i)\geq C(n, \nu_0, r_0, r)> 0 \label{3.3.6}
\end{align}

By (\ref{3.3.5}) and (\ref{3.3.6}), using compactness theorem for the Ricci flow in \cite{HamComp}, we can extract a subsequence Ricci flow solutions, which are convergent to a nonflat ancient smooth complete solution to  Ricci flow, which has max volume growth and $|Rm|\leq 4$. 

When $n=2$, it is obvious that $K_g^{\mathbb{C}}\geq 0$. We also have $K_{g}^{\mathbb{C}}\geq 0$ when $n= 3$, because any $3$-dim ancient solution has $Rm\geq 0$ by \cite{Chen}. When $n\geq 4$, $K_g^{\mathbb{C}}\geq 0$ is preserved. From Lemma $4.5$ in \cite{CW}, we know that for any complete nonflat ancient solution with bounded and nonnegative complex sectional curvature, its asymptotic volume ratio is $0$. That is contradiction with the max volume growth property we get in (\ref{volume lower bound}). Therefore Claim \ref{claim 3.3} is proved.
}
\qed

Without loss of generality, we can further assume $T_0\leq 1$. Now choose $r= 4r_0$ and apply Theorem \ref{thm chen}, we have $|Rm(x, t)|\leq C(n, \nu_0, r_0)$ on $B_{t}(x_0, \frac{3}{4}r_0)\times \in [0, T_0]$.
 
Consider the evolution equation of $u= \varphi(\frac{4d_t(x_0, x)}{r_0}- 1)|Rm(x, t)|^2$, where $\varphi$ is the smooth nonnegative  decreasing function chosen as before. 
\begin{equation}\nonumber
{\left.
\begin{array}{rl}
\square u &\leq \Big(\frac{4\varphi '}{r_0}\square d_{t}(x_0, x)- \frac{16}{r_0^2}\varphi ''\Big)|Rm|^2- 2\varphi |\nabla Rm|^2 \\
& -2\nabla \varphi \nabla |Rm|^2 + C_1\varphi |Rm|^3
\end{array} \right.
}
\end{equation}
where $C_1$ is some constant depending only on $n$. 

Let $u_{max}(t)= u(x(t), t)= \max_{x\in M^3}u(x, t)$, then by the maximum principle and the properties of $\varphi$, and using Lemma $8.3$ in \cite{Pere}, at $(x(t), t)$ we have
\begin{equation}\nonumber
{\left.
\begin{array}{rl}
\frac{d^{-}}{dt}u_{max}&\leq C(r_0)|Rm|^2+ C|Rm|u_{max}\\
&\leq C(\nu_0, r_0)+ C(\nu_0, r_0)u_{max}= C+ Cu_{max}
\end{array} \right.
}
\end{equation}
where $C= C(n, \nu_0, r_0)$ only depends on $n$, $\nu_0$ and $r_0$. Integrating this inequality, 
\begin{equation}\nonumber
{u_{max}(t)\leq [u_{max}(0)+ 1]e^{Ct}- 1\leq [r_0^{-4}+ 1]e^{Ct}- 1
}
\end{equation}

Hence there exists $C(n, \nu_0, r_0)> 0$ such that $u_{max}(t)< 4r_0^{-4}$ when $0\leq t\leq \min\{ C(n, \nu_0, r_0), T_0\}$. Then we get $|Rm(x, t)|< 2r_0^{-2}$ on $B_{t}(x_0, \frac{r_0}{2})$ when $0\leq t\leq \min\{C(n, \nu_0, r_0), T_0\}$. Recall that there is $(x_1, t_1)$ such that $t_1\leq T_0$, $x_1\in B_{t_1}(x_0, \frac{r_0}{2})$, $|Rm|(x_1, t_1)= 2r_0^{-2}$. Therefore $T_0\geq \min\{T, C(n, \nu_0, r_0) \}$. 
}
\qed

\section{Existence time estimate of Ricci flow}
Now we are ready to prove our main theorem.

{\it \textbf{Proof of Theorem \ref{thm 1.1}}:}~
{If $(M^n, g)$ is flat, We get our conclusion trivially. In the following of the proof, we assume $(M^n, g)$ is not flat. We have the following claim:
\begin{claim}\label{claim 3.1}
{$|Rm(x, t)|$ on $(M^n, g(t))$ blow up at the same time during the Ricci flow.
}
\end{claim}
\pf
{When $n= 2$, from Gauss-Bonnet Theorem we know that $M^2$ is diffeomorphic to $\mathbb{S}^2$. From \cite{Chow}, $(M^2, g)$ will be deformed to $\mathbb{S}^2$ with the canonical symmetric metric under the Ricci flow, the conclusion follows. If $n= 3$, from the main results in \cite{Ham} for positive Ricci curvature case and Sections $8$ and $9$ in \cite{Ham4} for the borderline case, $(M^3, g)$ will be deformed to locally symmetric metric space under the Ricci flow, and we get the conclusion. When $n\geq 4$, from results in \cite{Brendle} for positive case and a general maximum principle in \cite{BS2} applied for the borderline case, $(M^n, g)$ will be deformed to locally symmetric metric space under the Ricci flow, then the conclusion is obtained. From all the above, the claim is proved.
}
\qed

Assume on $(M^n, g)$, the Ricci flow's maximal existence time interval is $[0, T_1)$, where $0< T_1< \infty$. We define $T_0$ as the following:
\begin{align}
T_0\vcentcolon= \max{ \Big\{t\Big| \ t< T_1,\ |Rm|(x, s)\leq 2r_0^{-2},\ when \ x\in B_{t}\Big(x_0, \frac{1}{2}\epsilon_0 r_0\Big),\ and \ s\in [0, t],\  \Big\} } \nonumber 
\end{align}

From the above claim, we get $T_0< T_1$, hence the Ricci flow has the smooth solution on $M^n\times [0, T_0]$. 

If $T_0< C(r_0, n)$, where $C(r_0, n)$ is from Proposition \ref{prop 2.2}, then there exists $\epsilon> 0$ such that $T_0+ \epsilon \leq C(\nu_0, r_0)$ and the Ricci flow has smooth solution on $M^n\times [0, T_0+ \epsilon]$. From Proposition \ref{prop 2.2}, we get $|Rm|(x, t)\leq 2r_0^{-2}$ for $x\in B_{t}\Big(x_0, \frac{1}{2}\epsilon_0 r_0\Big)$ and $t\in [0, T_0+ \epsilon]$. It is contradiction with the definition of $T_0$. 

Hence $T_0\geq C(r_0, n)$, we get our conclusion.  
}
\qed

{\it \textbf{Proof of Theorem \ref{thm 1.2}}:}~
{Recall that $(M^n, g(x, 0))\times \mathbb{R}$ has nonnegative isotropic curvature is implied by nonnegative complex sectional curvature when $n\geq 4$ and is equivalent to $Rc\geq 0$ when $n= 3$. Hence we still have Claim \ref{claim 3.1} under the assumption of Theorem \ref{thm 1.2}. The rest proof of Theorem \ref{thm 1.2} is similar to the above proof except that we use Proposition \ref{prop 3.2} instead of Proposition \ref{prop 2.2}.
}
\qed

\begin{cor}\label{cor 4.2}
{Let $(M^3, g)$ be a complete noncompact $3$-dimensional manifold with $Rc\geq 0$. Assume there exists an exhaustion $\{\Omega_i\}_{i= 1}^{\infty}$ of $(M^3, g)$ with the property $\Omega_i\subset\subset \Omega_{i+ 1}$ for each $i$, and for each $\Omega_i$ there exists $(M_i, g^{(i)})$ with $Rc(g^{(i)})\geq 0$ such that $\Omega_i$ can be isometrically embedded into $(M_i, g^{(i)})$. Then there exists a constant $T> 0$ such that the Ricci flow has a complete solution on $M^3\times [0, T]$ with $g(0)= g$.
}
\end{cor}

\pf
{The proof of the corollary is similar to the proof of Theorem $4.7$ in \cite{CW}, except that we use Theorem \ref{thm 1.2} here instead of Propositions $4.3$ and $4.6$ there. Note in this case, the estimate (\ref{1.2.1}) is similar to Proposition $4.6$ there.
}
\qed

\begin{example}\label{Exam 1}
{From Proposition $6.7$ in \cite{CW}, if $(M^3, g)$ is a $3$-dim open manifold with $Rm\geq 0$, then there exists an exhaustion $\{\Omega_i\}_{i= 1}^{\infty}$ of $(M^3, g)$ with the property $\Omega_i\subset\subset \Omega_{i+ 1}$ for each $i$. And for each $\Omega_i$, there exists $(M_i, g^{(i)})$ with $Rm(g^{(i)})\geq 0$ such that $\Omega_i$ can be isometrically embedded into $(M_i, g^{(i)})$. By the Corollary above, the Ricci flow has a complete solution on $M^3\times [0, T]$ for some $T> 0$, this is the $3$-dim case of Theorem $1$ in \cite{CW}. 
}
\end{example}

\appendix
\section{Perelman's Lemma}\label{App 1}
The following lemma is originally due to G. Perelman (see section $10$ in \cite{Pere}). Following the updated version in \cite{RF3} (Chapters $21$ and $22$ there), we do some modifications to make it suitable for our use in the proof of Claim \ref{claim 3.3}. Although the modification is slight, for completeness, we provide the details for reader's convenience. A similar version of Lemma \ref{lem 3.1} and Claim \ref{claim 2.1} was proved independently in \cite{Simon3} (Theorem $2.1$ there), also see Lemmas $3.1.1$, $32.1$ in \cite{KL}. 

\begin{lemma}\label{lem 3.1}
{$(M^n, g(t))$ is a Ricci flow complete solution on $[0, T]$. Let $B>0$, $r> 0$ be some fixed constants, $\hat{x}$ is some point in $B_{g(\hat{t})}(x_0, \frac{r}{4})$ for some $\hat{t}\in (0, T]$, and satisfies $|Rm(\hat{x}, \hat{t})|_{g(\hat{t})}\geq B+ \frac{B}{\hat{t}}$. For any constant $A$ such that $1\leq A\leq B$, we can choose $(\bar{x}, \bar{t})$ with $\bar{Q}= |Rm(\bar{x}, \bar{t})|_{g(\bar{t})}\geq B\bar{t}^{-1}$ such that $\bar{x}\in B_{g(\bar{t})}(x_0, 2AB^{-\frac{1}{2}}+ \frac{r}{4})$ and 
\begin{align}
|Rm(x, t)|_{g(t)}\leq 4\bar{Q}\ , \quad \quad  (x, t)\in B_{g(\bar{t})}\Big(\bar{x}, \frac{1}{10}A\bar{Q}^{-\frac{1}{2}}\Big)\times \Big[\bar{t}- \frac{2}{25n}A\bar{Q}^{-1}, \bar{t}\Big] \nonumber
\end{align}

}
\end{lemma}

\pf
{Step (I). We firstly prove the following claim.

\begin{claim}\label{claim 2.1}
{We can find $(\check{x}, \check{t})$, $|Rm|(\check{x}, \check{t})\geq B\check{t}^{-1}$ and $\check{t}\in (0, T]$, $d_{g(\check{t})}(\check{x}, x_0)\leq 2AB^{-\frac{1}{2}}+ \frac{r}{4}$ such that
\begin{equation}\label{3.1.1}
{|Rm|(x, t)\leq 4|Rm|(\check{x}, \check{t})
}
\end{equation}
whenever
\begin{equation}\label{3.1.2}
{|Rm|(x, t)\geq Bt^{-1}, \  0< t\leq \check{t}, \  d_{g(t)}(x, x_0)\leq d_{g(\check{t})}(\check{x}, x_0)+ A|Rm|^{-\frac{1}{2}}(\check{x}, \check{t})
}
\end{equation}
For (\ref{3.1.1}) or (\ref{3.1.2}), we say $(x, t)$ satisfies (\ref{3.1.1}) or (\ref{3.1.2}) for $(\check{x}, \check{t})$.
}
\end{claim}

{\it \textbf{Proof of Claim \ref{claim 2.1}}:}~
{We use the notation $(x_1, t_1)$ for $(\hat{x}, \hat{t})$. If we choose $(\check{x}, \check{t})$ as $(x_1, t_1)$, and any $(x, t)$ satisfying (\ref{3.1.2}) will satisfy (\ref{3.1.1}) for $(x_1, t_1)$, then we are done. Otherwise, we can find $(x_2, t_2)$ satisfies (\ref{3.1.2}) but not (\ref{3.1.1}) for $(x_1, t_1)$. Then 
\begin{equation}\nonumber
{|Rm|(x_2, t_2)> 4|Rm|(x_1, t_1)\geq 4B, \quad |Rm(x_2, t_2)|\geq Bt_2^{-1}
}
\end{equation}
and 
\begin{equation}\nonumber
{d_{g(t_2)}(x_2, x_0)\leq d_{g(t_1)}(x_1, x_0)+ AB^{-\frac{1}{2}}\leq \frac{r}{4}+ AB^{-\frac{1}{2}}
}
\end{equation}
By induction, for any positive integer $k$, we can find $(x_k, t_k)$ such that
\begin{equation}\label{3.1.5}
{|Rm|(x_k, t_k)\geq 4^{k- 1}B\ , \quad \quad d_{g(t_k)}(x_k, x_0)\leq \frac{r}{4}+ \sum_{i= 2}^{k} A\cdot |Rm|^{-\frac{1}{2}}(x_{i- 1}, t_{i- 1})
}
\end{equation}
From (\ref{3.1.5}), we get 
\begin{align}
d_{g(t_k)}(x_k, x_0)&\leq \frac{r}{4}+ \sum_{i= 2}^{k} A\cdot |Rm|^{-\frac{1}{2}}(x_{i- 1}, t_{i- 1})\leq \frac{r}{4}+ \sum_{i= 2}^k A\cdot \big(4^{i- 2} B\big)^{-\frac{1}{2}} \nonumber \\
&\leq \frac{r}{4}+ 2AB^{-\frac{1}{2}} \label{domain bound}
\end{align}

Because $B_{g(t)}(x_0, \frac{r}{4}+ 2AB^{-\frac{1}{2}})\times [0, T]$ is compact, $|Rm|$ is bounded on it. From (\ref{domain bound}), $|Rm|(x_k, t_k)$ are all bounded. By the first inequality of (\ref{3.1.5}), the above induction stops at finite steps. Hence we find $(\check{x}, \check{t})$ such that $d_{g(\check{t})}(\check{x}, x_0)\leq 2AB^{-\frac{1}{2}}+ \frac{r}{4}$ and $|Rm|(\check{x}, \check{t})\geq B\check{t}^{-1}$. Furthermore, if $(x, t)$ satisfies (\ref{3.1.2}) for $(\check{x}, \check{t})$, then $(x, t)$ satisfies (\ref{3.1.1}) for $(\check{x}, \check{t})$. The claim is proved.
}
\qed

Step (II). Let $(\bar{x}, \bar{t})$ to be $(\check{x}, \check{t})$ found in the above claim. We prove the theorem by contradiction. Assume there exists some point $(\tilde{x}, \tilde{t})$, such that
\begin{align}
(\tilde{x}, \tilde{t})\in B_{g(\bar{t})}\Big(\bar{x}, \frac{1}{10}A\bar{Q}^{-\frac{1}{2}}\Big)\times \Big[\bar{t}- \frac{2}{25n}A\bar{Q}^{-1}, \bar{t}\Big] \nonumber
\end{align}
\begin{align}
|Rm|(\tilde{x}, \tilde{t})> 4|Rm|(\bar{x}, \bar{t}) \label{ineq 0}
\end{align}

Using the assumptions $\bar{Q}\geq B\bar{t}^{-1}$ and $A\leq B$, we get 
\begin{align}
\bar{t}- \frac{2}{25n}A\bar{Q}^{-1}\geq \bar{t}- \frac{2}{25n}\bar{t}> \frac{1}{4}\bar{t} \label{ineq 1}
\end{align}

From (\ref{ineq 0}), (\ref{ineq 1}) and $\tilde{t}\geq \bar{t}- \frac{2}{25n}A\bar{Q}^{-1}$, we get $|Rm|(\tilde{x}, \tilde{t})> 4B\bar{t}^{-1}> B\tilde{t}^{-1}$. 

From triangle inequality, $\bar{B}_{g(\bar{t})}\Big(\bar{x}, \frac{1}{10}A\bar{Q}^{-\frac{1}{2}}\Big)\subset B_{g(\bar{t})}\Big(x_0, d_{g(\bar{t})}(\bar{x}, x_0)+ \frac{9}{10}A\bar{Q}^{-\frac{1}{2}}\Big)$. Suppose that there exists $t_1\geq \bar{t}- \frac{2}{25n}A\bar{Q}^{-1}$ which is the first time , going backwards in time from $\bar{t}$, such that the boundary of $\bar{B}_{g(t_1)}\Big(x_0, d_{g(\bar{t})}(\bar{x}, x_0)+ \frac{9}{10}A\bar{Q}^{-\frac{1}{2}}\Big)$ intersects $\bar{B}_{g(\bar{t})}\Big(\bar{x}, \frac{1}{10}A\bar{Q}^{-\frac{1}{2}}\Big)$. Let $x_{*}$ be such a point of intersection, then
\begin{equation}\label{3.1.8}
{d_{g(t_1)}(x_{*}, x_0)= d_{g(\bar{t})}(\bar{x}, x_0)+ \frac{9}{10}A\bar{Q}^{-\frac{1}{2}}
}
\end{equation}
By $x_{*}\in \bar{B}_{g(\bar{t})}(\bar{x}, \frac{1}{10}A\bar{Q}^{-\frac{1}{1}})$, we get
\begin{equation}
{d_{g(\bar{t})}(x_{*}, x_0)\leq d_{g(\bar{t})}(\bar{x}, x_0)+ \frac{1}{10}A\bar{Q}^{-\frac{1}{2}} \label{used later}
}
\end{equation} 
From $Rc\geq 0$ and the Ricci flow equation, $d_{g(t)}(x_{*}, x_0)\leq d_{g(t_1)}(x_{*}, x_0)$ for $t\in [t_1, \bar{t}]$. Combining with (\ref{3.1.8}), we get
\begin{align}
d_{g(t)}(x_{*}, x_0)\leq d_{g(\bar{t})}(\bar{x}, x_0)+ \frac{9}{10}A\bar{Q}^{-\frac{1}{2}} \ , \quad \quad t\in [t_1, \bar{t}] \label{C_1}
\end{align}
Then from (\ref{C_1}) and triangle inequality,
\begin{equation}\label{C_2}
{\bar{B}_{g(t)}\Big(x_{*}, \frac{1}{10}A\bar{Q}^{-\frac{1}{2}}\Big)\subset \bar{B}_{g(t)}\Big(x_0, d_{g(\bar{t})}(\bar{x}, x_0)+ A\bar{Q}^{-\frac{1}{2}}\Big)\ , \quad \quad t\in [t_1, \bar{t}]
}
\end{equation}
When $x\in \bar{B}_{g(t)}\Big(x_{0}, d_{g(\bar{t})}(\bar{x}, x_0)+ A\bar{Q}^{-\frac{1}{2}}\Big)$, $t\in \Big[\bar{t}- \frac{2}{25n}A\bar{Q}^{-1}, \bar{t}\Big]$, if 
$|Rm|(x, t)< Bt^{-1}$, then using (\ref{ineq 1}),
\begin{align}
|Rm|(x, t)< Bt^{-1}\leq 4B\bar{t}^{-1}= 4\bar{Q} \label{C_1_1}
\end{align}
If $|Rm|(x, t)\geq Bt^{-1}$, by Claim \ref{claim 2.1} we get $|Rm|(x, t)\leq 4\bar{Q}$ again, combining with (\ref{C_1_1}) we have
\begin{align}
|Rm|(x, t)\leq 4\bar{Q}\ , \  x\in \bar{B}_{g(t)}\Big(x_{0}, d_{g(\bar{t})}(\bar{x}, x_0)+ A\bar{Q}^{-\frac{1}{2}}\Big)\ , \ t\in \Big[\bar{t}- \frac{2}{25n}A\bar{Q}^{-1}, \bar{t}\Big] \label{C_3}
\end{align}
From (\ref{C_2}) and (\ref{C_3}),
\begin{align}
|Rc(x, t)|\leq 4(n- 1)\bar{Q}\ , \quad  x\in \Big\{\bar{B}_{g(t)}(x_{*}, \frac{1}{10}A\bar{Q}^{-\frac{1}{2}})\cup \bar{B}_{g(t)}(x_0, \frac{1}{10}A\bar{Q}^{-\frac{1}{2}})\Big\} \ , \ t\in [t_1, \bar{t}] \nonumber
\end{align} 
Hence we can apply Lemma $8.3\ (b)$ in \cite{Pere} (where we choose $K= 4\bar{Q}$, $r_0= \bar{Q}^{-\frac{1}{2}}$), we have 
\begin{equation}
{\frac{d}{dt}d_{g(t)}(x_{*}, x_0)\geq -10(n- 1)\bar{Q}^{\frac{1}{2}} \label{ODE}
}
\end{equation}
Integrate (\ref{ODE}), using (\ref{used later}), we get
\begin{equation}\nonumber
{\left.
\begin{array}{rl}
d_{g(t_1)}(x_{*}, x_0)&\leq d_{g(\bar{t})}(x_{*}, x_0)+ 10(n- 1)\bar{Q}^{\frac{1}{2}}\cdot\frac{2}{25n}A\bar{Q}^{-1}  \\
& < d_{g(\bar{t})}(\bar{x}, x_0)+ \frac{9}{10}A\bar{Q}^{-\frac{1}{2}}
\end{array} \right.
}
\end{equation}
This is the contradiction with (\ref{3.1.8}). Hence for $t\in \Big[\bar{t}- \frac{2}{25n}A\bar{Q}^{-\frac{1}{2}}, \bar{t}\Big]$,
\begin{align}
\bar{B}_{g(\bar{t})}(\bar{x}, \frac{1}{10}A\bar{Q}^{-\frac{1}{2}})\subset \bar{B}_{g(t)}\Big(x_0, d_{g(\bar{t})}(\bar{x}, x_0)+ \frac{9}{10}A\bar{Q}^{-\frac{1}{2}}\Big)\label{last one}
\end{align}
By $\tilde{x}\in \bar{B}_{g(\bar{t})}(\bar{x}, \frac{1}{10}A\bar{Q}^{-\frac{1}{2}})$, choose $t= \tilde{t}$ in (\ref{last one}), we get
\begin{equation}\nonumber
{d_{g(\tilde{t})}(\tilde{x}, x_0)\leq d_{g(\bar{t})}(\bar{x}, x_0)+ \frac{9}{10}A\bar{Q}^{-\frac{1}{2}}
}
\end{equation}
By Claim \ref{claim 2.1}, we get $|Rm|(\tilde{x}, \tilde{t})\leq 4|Rm|(\bar{x}, \bar{t})$. That is contradiction.
}
\qed

\section*{Acknowledgments}
The author would like to thank Bing-Long Chen and Miles Simon for their comments. He is grateful to Gang Liu for pointing out the possible connection between Perelman's pseudolocality and the results in the earlier version of this note. He is also indebted to Simon Brendle for discussion and comments, Yuan Yuan, Zhiqin Lu and Jeff Streets for suggestions. Finally, he is deeply grateful to the referee for careful reading and for many helpful suggestions, which greatly improved the readability of the note.


 

\end{document}